\documentclass[doublecol]{epl2}
\usepackage[utf8]{inputenc}
\usepackage{graphicx}
\usepackage{xcolor}
\usepackage{hyperref}
\usepackage{amsthm,amsfonts,amsmath,amssymb}

\newtheorem{algorithm}{Algorithm}

\theoremstyle{definition}

\newcommand{\prob}{\mathbb{P}}
\newcommand{\Prob}[1]{\prob\left(#1\right)}

\title{Switch chain mixing times through triangle counts}
\author{Tom Bannink\inst{1} \and Remco van der Hofstad\inst{2} \and Clara Stegehuis\inst{2}}
\shortauthor{Tom Bannink, Remco van der Hofstad and Clara Stegehuis}

\institute{                    
  \inst{1} QuSoft and CWI, Amsterdam, The Netherlands\\
  \inst{2} Eindhoven University of Technology, Department of Mathematics and Computer Science, Eindhoven, The Netherlands
}

\abstract{
    Sampling uniform simple graphs with power-law degree distributions with degree exponent $\tau\in(2,3)$ is a non-trivial problem. We propose a method to sample uniform simple graphs that uses a constrained version of the configuration model together with a Markov Chain switching method. We test the convergence of this algorithm numerically in the context of the presence of small subgraphs. We then compare the number of triangles in uniform random graphs with the number of triangles in the erased configuration model. Using simulations and heuristic arguments, we conjecture that the number of triangles in the erased configuration model is larger than the number of triangles in the uniform random graph, provided that the graph is sufficiently large. 
}

\begin{document}
\maketitle

\section{Introduction}
Many real-world networks have been found to have a power-law degree distribution with exponent $\tau\in(2,3)$~\cite{albert1999,faloutsos1999,vazquez2002}.
A uniform random graph with prescribed degrees serves as a null model for real-world networks, and has attracted enormous attention in network physics~\cite{newman2001,roberts2000,bender1978}. The uniform random graph (URG) is a uniform sample from the ensemble of all possible graphs with the prescribed degree sequence. 
The configuration model is used frequently to generate URGs~\cite{bollobas1980}. The configuration model starts with $n$ vertices and a degree sequence $(d_i)_{i=1,\dots, n}$ such that the sum of the degrees is even. All vertices $i$ start with $d_i$ half-edges, where $d_i$ is the degree of vertex $i$. Then, these half-edges are paired one by one, uniformly at random. This creates a random graph with the desired degree distribution. When the configuration model results in a simple graph, this is a uniform sample of all simple graphs with that degree sequence. As long as the degree exponent $\tau$ satisfies $\tau>3$, the probability that the configuration model creates a simple graph is strictly positive. Thus, in this regime, the configuration model can be effectively repeated until it results in a simple graph. When $\tau<3$, the probability that the configuration model results in a simple graph vanishes instead. Thus, the configuration model cannot be used to generate uniform simple graphs for $\tau\in(2,3)$.  

Several models exist to generate graphs with \emph{approximately} the desired degree sequence~\cite{britton2006}. One such model is the erased configuration model (ECM), where after the construction of the configuration model, all self-loops and multiple edges are removed. Another option is to use models with soft constraints on the degrees, such as hidden-variable models~\cite{chung2002, boguna2003}. Other methods to sample uniform graphs are based on maximizing entropy~\cite{squartini2015}. 

A method to sample random graphs with \emph{exactly} the desired degree sequence is to use Markov Chains. These methods start with an initial graph with the desired degree sequence. Then, at every time step, some edges of the graph are rewired in such a way that the stationary distribution is uniform~\cite{milo2003, coolen2009, greenhill2017, carstens2017}. When the number of rewirings (or \emph{switches}) tends to infinity, the result is a uniformly sampled random graph from all simple graphs with the same degree sequence.
These Markov Chain methods can be adapted to generate directed graphs~\cite{artzy2005,roberts2012}, connected graphs~\cite{mihail2003} or graphs with fixed degree-degree correlations~\cite{roberts2012}. 
 
There are several ways of creating an initial graph with the desired degree sequence, one of which is by using the classical Havel-Hakimi algorithm. We introduce a new algorithm that uses a constrained version of the configuration model where self-loops and multiple edges are avoided. We experimentally study the effect of the choice of initial graph on the number of switches needed to reach equilibrium. We show that the new algorithm does not let the Markov Chain produce uniform random graphs any faster than the Havel-Hakimi algorithm.

We analyze the influence of the starting configuration in the context of the presence of triangles, similar to~\cite{milo2003}. Triangles are the smallest nontrivial subgraphs of networks, and indicate the presence of communities or hierarchies~\cite{colomer2013,ravasz2003} or geometry in networks~\cite{krioukov2016} and influence the behaviour of spreading processes of networks~\cite{serrano2006}. We therefore experimentally study the number of Markov Chain switches required until the density of the number of triangles reaches equilibrium from different starting states.

In the ECM with $\tau\in(2,3)$, the number of triangles scales as $n^{\frac{3}{2}(3-\tau)}$~\cite{hofstad2017d}. We numerically investigate the scaling of the number of triangles in URG, using the switch chain. We find that finite-size effects play a role at $n=10.000$, and the data suggests that the scaling present in the ECM is the same in URG, but with different multiplicative constants.

\section{Switch chain}
We now explain the Markov Chain switching method we study in more detail. Its state space is the space of simple graphs with the desired degree distribution. At every time step $t$, it selects two vertex-distinct edges of the graph uniformly at random, say $(u_1,v_1)$ and $(u_2,v_2)$. These edges are replaced by $(u_1,v_2)$ and $(u_2,v_1)$ if this results in a simple graph, and otherwise the switch is rejected so that the graph remains the same. In both cases we set $t=t+1$. 
Setting $t=t+1$ also when a switch is rejected is crucial: if we do not increase the time after a rejected switch, the stationary state of the switch chain may not be uniform~\cite{carstens2017}. On the other hand, when we do increase the time step after a rejected switch, the stationary state of the switch chain is the uniform distribution. See also section A of the supplementary material.
When the degree sequence behaves like a power law, bounds on the mixing time of the switch chain are large (for example $n^9$ in~\cite{greenhill2017}) when the degree exponent $\tau>3$, and unknown in the case where $\tau<3$. Experimental results suggest that the mixing time is much smaller than the bounds that have been proven~\cite{rechner2016,ray2014}.  

\begin{figure}
    \centering
    \includegraphics[width=0.45\textwidth]{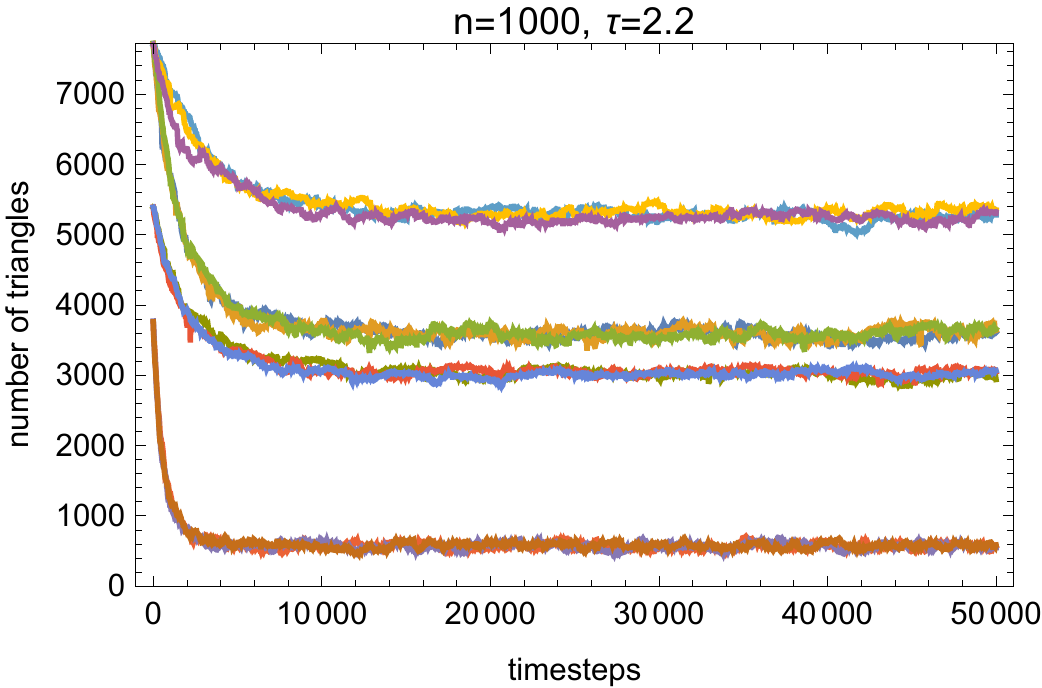} 
    \caption{\label{fig:timeevol} Time evolution of the number of triangles in the switch chain ($n=1000$, $\tau=2.2$). The horizontal axis shows steps of the Markov Chain and the vertical axis shows the number of triangles. Four degree sequences sampled from $D$ given by \eqref{eq:degreedistr} are shown. For each degree sequence, the Havel-Hakimi construction was performed and three runs of the switch chain with that same starting point are shown. The overlapping lines in the plot correspond to the different runs on the same starting graph.}
\end{figure}

\section{Degree distribution}
We study networks with power law degree distribution $D$ in the infinite variance regime, so that $D$ satisfies
\begin{equation}\label{eq:degreedistr}
    \Prob{D=k}\approx C k^{-\tau},
\end{equation}
for some constant $C$ and with $\tau\in(2,3)$ and $k$ large.
Several methods exist to sample a degree sequence from $D$. For example, it is possible to sample $n$ i.i.d.\ copies from~\eqref{eq:degreedistr}. Then, the variability in the degree sequence is the largest contributor to the variance of several network observables~\cite{ostilli2014}. This is also visible in Figure \ref{fig:timeevol}, where the number of triangles varies enormously between different sampled degree sequences from the same distribution. However, when we create a null model corresponding to a real-world network, the degree sequence of this null model is usually fixed, and there is no variability in the degree sequence. Therefore, we sample the power-law degree sequence by using the inverse of the distribution function $F_\tau$ of $D$.

Let us explain how this can be done. Define
\begin{equation}\label{eq:degreefixed}
    d^{(n,\tau)}_i=[1-F_\tau]^{-1}(i/n)
\end{equation}
for $i\in[n]$. Then, the empirical degree distribution converges in distribution to $F_\tau$ as $n\to\infty$. Note that the degrees $d^{(n,\tau)}_i$ are indeed deterministic. Then the only uncertainty in the resulting random graph is from the random connection of the edges, similar to what we encounter when creating a null model for a particular observation of a real-world network. 
We will refer to the degree sequence $(d^{(n,\tau)}_i)_{i\geq 1}$ given by \eqref{eq:degreefixed} as the \emph{canonical degree sequence} for a given $(n,\tau)$, see also Chapter 7 of~\cite{hofstad2009}. See section B of the supplementary material for more details on sampling from $D$ and the canonical degree sequence.

\section{Initial graphs for the switch chain}
The switch chain needs an initial graph to start switching. We investigate three different methods to obtain a simple graph with a given degree distribution, and see what the effect is of this initial graph on the performance of the switch chain.

\subsection{Constrained configuration model (CCMd)}
The constrained configuration model (CCMd) is defined as follows:
\begin{algorithm}[Constrained configuration model] ~\\
\textbf{Input:} A degree sequence $(d_1,\dots,d_n)$ with corresponding vertices $(v_1,\dots,v_n)$.\\
\textbf{Output:} A simple graph with degree sequence $(d_1,\dots,d_n)$ or \textsc{fail}.
\begin{enumerate}
    \item Let $V=\{v_1,\dots,v_n\}$ be the set of vertices. Set $W=\emptyset$. Equip vertex $i$ with $d_i$ half-edges for every $i$. 
    \item \textbf{while} there are half-edges \textbf{do}
    \begin{enumerate}
        \item Let $W=\{v\}$ where $v$ is the vertex with the highest amount of remaining half-edges incident to it.
        \item \textbf{while} $v$ has half-edges \textbf{and} $V\setminus W$ has incident half-edges \textbf{do}
        \begin{itemize}
            \item Pair a half-edge adjacent to $v$ to a uniformly chosen half-edge adjacent to $V\setminus W$. Denote the vertex to which $v$ is paired by $w$ and remove both half-edges.
            \item Set $W=W\cup \{w\}$.
        \end{itemize}
        \item If $v$ has unpaired half-edges, output \textsc{fail}.
    \end{enumerate}
\end{enumerate}
\end{algorithm}

Thus, the algorithm works as the configuration model, except that it keeps track of a list $W$ of `forbidden vertices' that guarantees that no self-loops or multiple edges are created. Note that this algorithm may fail and not produce a simple graph with the desired degree sequence. For example, the last vertex may have two unpaired half-edges incident to it. Then, the only way to finish the pairing is to create a self-loop, which we have forbidden. First choosing the vertex $v$ with the highest number of half-edges aims at avoiding the algorithm to fail: pairing the highest-degree vertices without conflicts is the most difficult. When we pair these vertices at the start of the algorithm, the probability that these are paired successfully is larger. Note that this algorithm does not create a uniformly sampled simple graph, as the regular configuration model would. 

\subsection{Constrained configuration model, updated (CCMdu)}
A variation on the constrained configuration model is the updated constrained configuration model (CCMdu). Where the constrained configuration model algorithm pairs all half-edges incident to the chosen vertex $v$ before proceeding to the next vertex $v'$, the updated constrained configuration model only does one pairing before replacing $v$ by the vertex with the highest amount of remaining half-edges:
\begin{algorithm}[Updated constrained configuration model] ~\\
\textbf{Input:} A degree sequence $(d_1,\dots,d_n)$ with corresponding vertices $(v_1,\dots,v_n)$.\\
\textbf{Output:} A simple graph with degree sequence $(d_1,\dots,d_n)$ or \textsc{fail}.
\begin{enumerate}
    \item Let $V=\{v_1,\dots,v_n\}$ be the set of vertices. Let $W_i=\{v_i\}$ for all $i\in 1,\dots,n$. Equip vertex $j$ with $d_j$ half-edges for all $j$. 
    \item \textbf{while} there are half-edges \textbf{do}
    \begin{enumerate}
        \item Let $v_i$ be the vertex with the highest amount of unpaired half-edges incident to it.
        \item \textbf{if} there are unpaired half-edges in $V\setminus W_i$ \textbf{then}\\
              Pair a half-edge adjacent to $v_i$ to a uniformly chosen half-edge adjacent to $V\setminus W_i$ and remove both half-edges. Denote the vertex to which $v_i$ is paired by $w$. Set $W_i=W_i\cup\{w\}$. \\
              \textbf{else} output \textsc{fail}
    \end{enumerate}
\end{enumerate}
\end{algorithm}
Just like the previous algorithm, this algorithm is not guaranteed to finish successfully, and does not create uniformly sampled graphs. 

\subsection{Havel-Hakimi - Erd\H{o}s-Gallai}
The Havel-Hakimi algorithm is a simple deterministic algorithm to create simple graphs. This algorithm sorts the degree sequence as $d_1\geq d_2\geq\dots\geq d_n$. Then it pairs vertex $v_1$, with the highest degree $d_1$, to $v_2,\dots ,v_{d_1+1}$. The degrees of these vertices are reduced by 1 (vertex $v_1$ is now done) and the procedure is repeated by re-sorting the degrees and pairing the new vertex with the highest degree. The Erd\H{o}s-Gallai theorem states that this algorithm always finishes in a simple graph with the desired degree sequence if such a graph exists.

Note that this graph is highly unlike a uniform sample of all graphs with the same degree sequence. The Havel-Hakimi construction creates many triangles and complete graphs of other sizes, see also section C of the supplementary material. In uniform samples, we expect to see fewer triangles and such larger complete graphs.

\section{(Empirical) mixing time}
As stated before, the switch chain mixing time for degree sequences sampled from a power-law degree distribution with $\tau\in(2,3)$ is unknown. Yet, we want to stop the switch chain at some point and get a sample graph from the uniform distribution. A common thing to do is computing graph properties like clustering coefficients, number of triangles, diameters and graph eigenvalues \cite{milo2003,ray2012}. The number of triangles in a network is an important observable, since it indicates the presence of communities or hierarchies~\cite{colomer2013,ravasz2003} or geometry in networks~\cite{krioukov2016} and it influences the behaviour of spreading processes of networks~\cite{serrano2006}. We therefore study the time evolution of the number of triangles in the switching process and stop when this quantity has sufficiently stabilised, similar to~\cite{milo2003}.


Figure \ref{fig:timeevol} shows the time evolution of the number of triangles for several samples of a degree distribution for $n=1000$.
First of all the figure shows that the number of triangles is highly dependent on the degree sequence, even if sampled from the same distribution. For runs of the switch chain with the same degree sequence, the number of triangles seems to evolve very similarly and in each run it takes about the same time for the number of triangles to stabilize.

\begin{figure}
    \centering
    \includegraphics[scale=1.0]{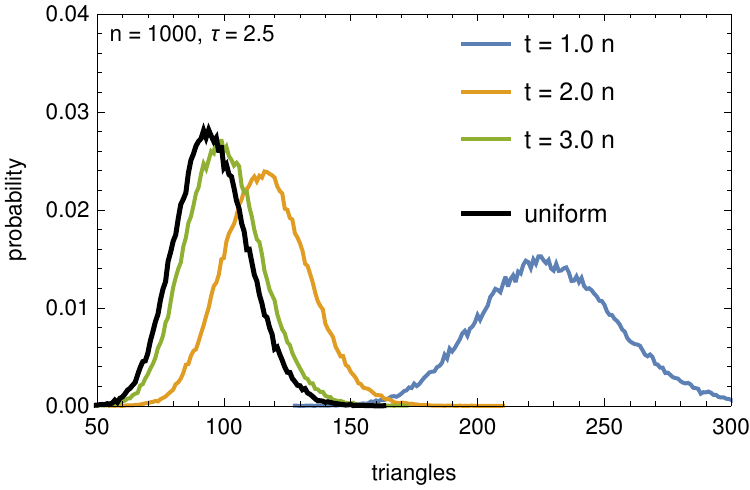}
    \caption{\label{fig:etmt} The distribution of triangles at different timesteps of the switch chain using the same (Havel-Hakimi) starting point. The degree sequence is the canonical degree sequence for $n=1000$ and $\tau=2.5$. The data is obtained by recording the number of triangles at $t=0.1n,\;t=0.2n, ...,\; t=20n$ and repeating this procedure 100.000 times. For the approximate uniform sample we record the number of triangles at $t=2000 n$.}
\end{figure}
To quantify the time it takes for the number of triangles to stabilize, we have computed an estimate of the distribution of the number of triangles at several timesteps for runs with the same starting point. One instance of this is shown in Figure \ref{fig:etmt}. For each of these distributions we have computed the total variation distance between it and the uniform distribution, and we consider the number of steps that it takes for the distance to become less than $0.1$.
Based on simulations up to $n=20.000$ we conclude that this empirical mixing time is at most $O(n\log^2 n)$ for constant $\tau\in(2,3)$.



\section{The number of triangles in uniform graphs}
\begin{figure}
    \centering
    \includegraphics[scale=0.8]{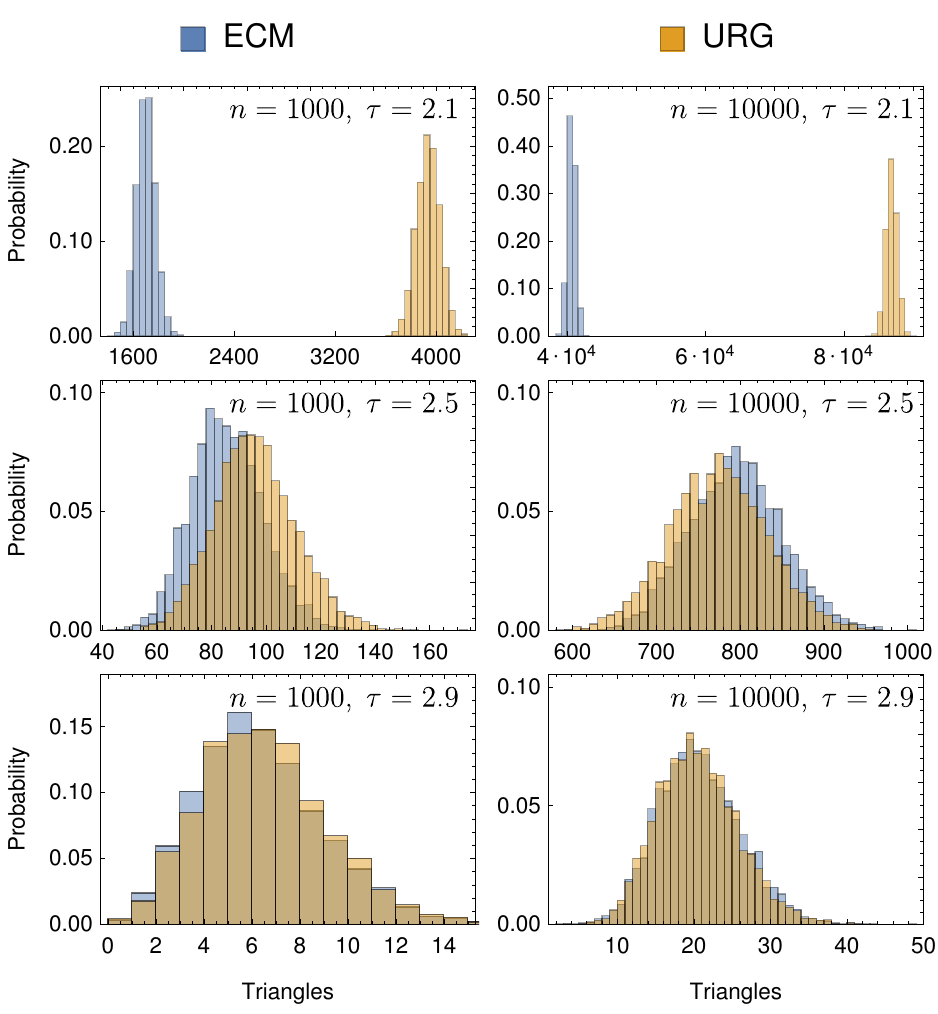}
    \caption{Distribution of the number of triangles in the ECM and URG. For every $n,\tau$, the canonical degree sequence was used and the ECM construction was performed 5000 times, shown together with 5000 samples from the uniform distribution.}
    \label{fig:ecmdistribution}
\end{figure}
In the ECM the number of triangles scales as $n^{\frac{3}{2}(3-\tau)}$~\cite{hofstad2017d}. To compare, Figure~\ref{fig:ecmdistribution} shows the distribution of the number of triangles for many values of $n$ and $\tau$, for both the ECM and the URG.
The figure shows that for $\tau$ close to 2, the ECM contains fewer triangles than the URG. However, for $n=10.000$ and $\tau=2.5$ the average number of triangles in the ECM is higher than in a uniform random graph. It is known that for $\tau>3$ these models are similar (see for example~\cite[Chapter 7]{hofstad2009}), so that for $\tau=2.9$, the difference between the number of triangles in the ECM and the URG is expected to be small.

We predict that, for $\tau\in(2,3)$, for $n$ very large, the number of triangles in the ECM is on average \emph{larger} than that in the URG, but this may only be visible for $n$ extremely large (see Figure~\ref{fig:ecmdistribution} for $\tau=2.5$). This prediction is quite counterintuitive, as the number of edges in the ECM is \emph{smaller} than that in the URG. This prediction is due the fact that the edge probabilities in the ECM are close to $1-\mathrm{e}^{-d_id_j/L_n}$ (see \cite{hofstad2005,hofstad2017c}), while in the URG it is $\frac{d_id_j}{L_n+d_id_j}$ (see \cite{mckay2010b}), which is a little larger. Here $d_i$ are the degrees and $L_n$ is the sum over all degrees. Further, the edges are close to being independent, so that the triangle counts are close to those in a hidden-variable model with the vertex weights being given by the degrees.
For $\tau>3$, this is worked out in detail in \cite[Chapter 7]{hofstad2009}.
In fact, the clustering coefficient for the hidden-variable model with connection probabilities $\frac{d_i d_j}{L_n+d_i d_j}$ is lower than the clustering coefficient for the hidden-variable model with connection probabilities $\mathrm{e}^{-d_i d_j / L_n}$~\cite{hofstad2017b}.
Intuitively, high degree vertices in uniform random graphs are forced to connect to low degree vertices, because otherwise the simplicity constraint on the graph would be violated. These low degree vertices barely participate in triangles. In the ECM, high degree vertices will be connected more frequently. Because high degree vertices participate in more triangles, this suggests that the ECM contains more triangles than a uniform random graph. However, we also see that this effect can only kick in for very high $n$, particularly when $\tau$ is close to 2 or 3. This explains why the effect is only visible for $n=10,000$ for $\tau=2.5,$ and not for the other two values. For $\tau=3$, the convergence is slow since the number of triangles grows like $n^{3(3-\tau)/2},$ and the exponent is small, while only the front factor is asymptotically different, which takes long to be seen. For $\tau$ close to 2, the convergence of the number of triangles to the large network limit is extremely slow in the hidden-variable model~\cite{hofstad2017b}, which may explain why the asymptotics also kicks in late for $\tau$ close to 2 in the URG and the ECM. Determining the size of $n$ when the URG starts having more triangles than the ECM, as a function of $\tau$, remains of substantial interest.

\begin{figure}
    \centering
    \includegraphics[scale=1.0]{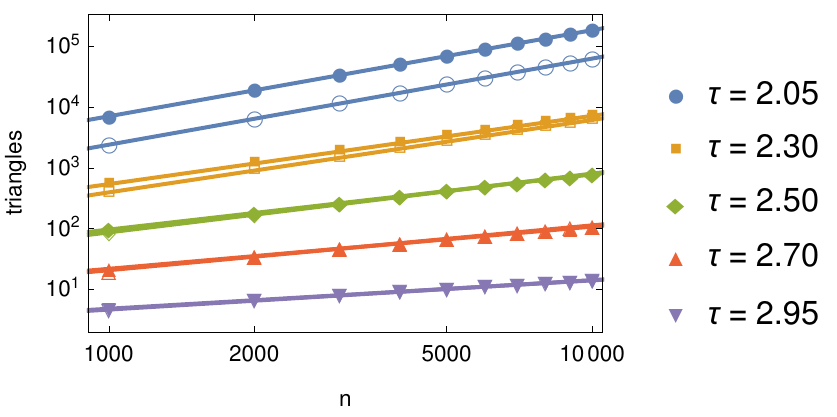}
    \caption{The datapoints show the average number of triangles in a log-log plot for several values of $n$ and $\tau$ for the canonical degree sequence. The solid markers correspond to the URG distribution, and the open markers correspond to the \emph{ECM}. The lines show a fit of the function $\log(\mathrm{triangles}) = a \log(n) + b$.}
    \label{fig:avgtrisloglog}
\end{figure}

\begin{figure}
    \centering
    \includegraphics[width=0.45\textwidth]{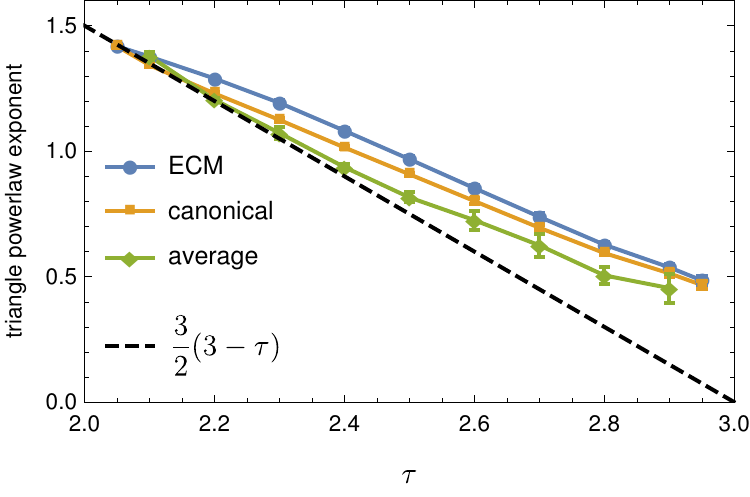}
    \caption{Values of the exponent in the triangle power law. The dashed line is $\frac{3}{2}(3-\tau)$ which is the theoretical exponent for the number of triangles in the ECM. The blue orange lines were obtained from fitting the data shown in Figure \ref{fig:avgtrisloglog}. The line labelled average is from a similar process but where the average was taken over 2000 sampled degree sequences instead of the canonical degree sequence. The error bars show the uncertainty of the fit parameters without taking into account the uncertainty of the data points themselves.}
    \label{fig:triangle_exponent}
\end{figure}

We considered the canonical degree sequence given by \eqref{eq:degreefixed} as well as $2000$ samples of valid degree sequences sampled from the distribution in \eqref{eq:degreedistr}. The results for the canonical degree sequence are shown in Fig. \ref{fig:avgtrisloglog}.
For $\tau=2.3$ the figure supports the conjecture as indeed the number of triangles in the ECM overtakes that of the URG as $n$ becomes large.
For each value of $\tau$ we have fitted the function $\log(\mathrm{triangles})=a\cdot\log(n)+b$ to this data to obtain Figure \ref{fig:triangle_exponent}.
As the figure shows, the data points lie above the line $\frac{3}{2}(3-\tau)$. Since $\frac{3}{2}(3-\tau)$ is known to be the correct exponent for the ECM, this also supports the prediction that the asymptotics only kick in at very large $n$.

\section{Switch chain as a proof method}
The switch chain has also been used as a combinatorial method for counting triangles~\cite{mckay2004,mckay2010,gao2012} in uniform random graphs. In these works, variations of the switch chain are studied where different edge rewiring rules are used but the Markov Chain still converges to the desired uniform distribution over the graphs. The idea is to count the number of triangles that a move of the Markov Chain can create or destroy when a switch is performed on certain vertices. Such proofs usually~\cite{mckay2004,gao2012} assume that such a move only creates or destroys at most one triangle with high probability.
\begin{figure}
    \centering
    \includegraphics[width=0.49\textwidth]{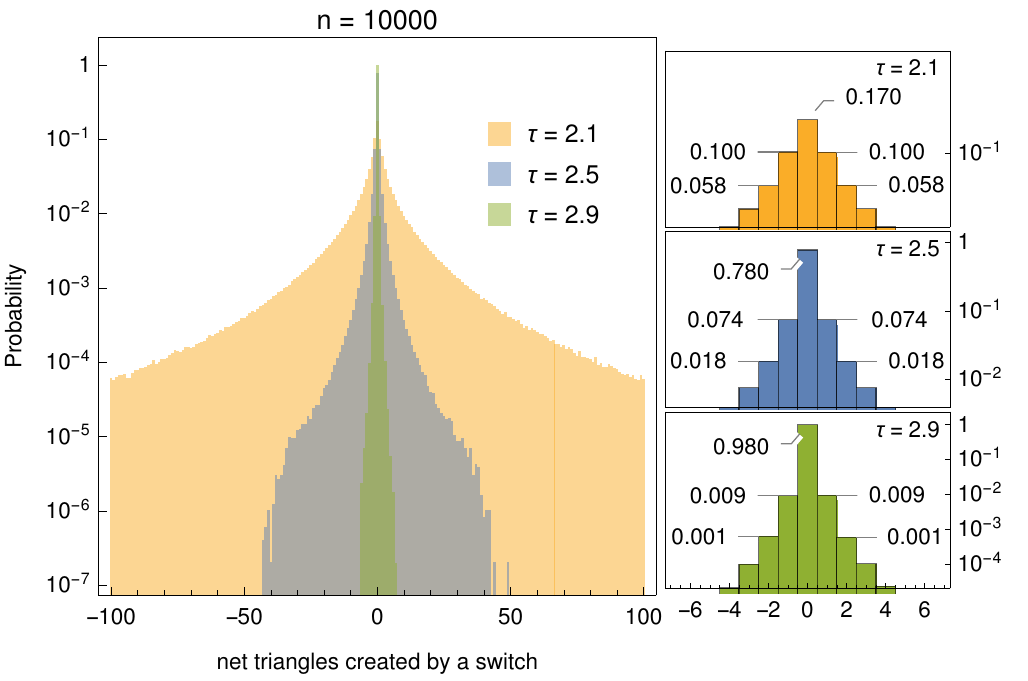}
    \caption{\label{fig:creationfreq} Triangle creation and destruction frequencies in equilibrium for $n=10000$ and three values of $\tau$ using the canonical degree sequence. Rejected switches are not counted here. The plots on the right are zoomed-in versions of the same plot. Note that it is possible that a move creates $l$ triangles and destroys $m$ other triangles which shows up as a net creation of $l-m$ triangles in these plots.}
\end{figure}
These proofs do not directly apply in the regime $\tau\in(2,3)$, so here we investigate this assumption numerically. Figure \ref{fig:creationfreq} shows histograms of the number of triangles that were created or destroyed by the switch chain moves in equilibrium. The plots show that the probability of creating or destroying $k>0$ triangles becomes lower as $\tau$ increases. It is important to keep in mind that these plots only show the net number of created triangles, so the proportion of moves that create or destroy any triangle might be considerably higher. We see that the probability that 2 or more triangles are created or destroyed can be large, especially when $\tau$ approaches 2. This suggests that these types of switch chain proofs cannot be used to count triangles for $\tau\in(2,3)$.

\section{Constrained configuration Model}
In this section we discuss the results for the two variants of the constrained configuration Model (CCM). Unlike the Havel-Hakimi construction, the CCM construction has a non-zero probability of failing and thereby not producing a simple graph with the desired degree sequence. On the other hand, the always-successful Havel-Hakimi starting state is far from uniform and it might be that the CCM construction provides a starting state that requires less switch chain moves to get good samples. We wish to investigate this and see whether the overhead of the CCM constructions is worth the (computational) time.
In this section we use the abbreviations CCMd and CCMdu to distinguish the two constructions introduced earlier.

\subsection{Construction success rate}
We looked at the construction success rates for CCMd and CCMdu, for graphical degree sequences. This means we only look at degree sequences for which a simple graph exists, checked using the Erd\H{o}s-Gallai theorem. It turns out that the construction success rate for CCMdu was lower than that of CCMd. For $n=1000$ and several values of $\tau$, we sampled 200 graphical degree sequences from the distribution given by \eqref{eq:degreedistr} and did 1000 CCMd construction attempts per sequence. We only found a single degree sequence with any failed attempts (4 out of 1000). All other degree sequences had 0 failed attempts and always successfully produced simple graphs.
For CCMdu, however, this is not the case. Although for most degree sequences it had a high success rate, there were degree sequences for $\tau$ close to 2, for which over 95\% of the attempts failed. This is further discussed in section E of the supplementary material.
Comparing the CCMd and CCMdu algorithms, we can conclude that finishing all pairings of a single vertex (CCMd) yields higher success rates. In a way CCMd is more similar to the Havel-Hakimi construction than CCMdu, because the Havel-Hakimi construction also finishes one vertex completely before moving on.

\subsection{Number of triangles}
\begin{figure}
    \centering
    \includegraphics[scale=1.0]{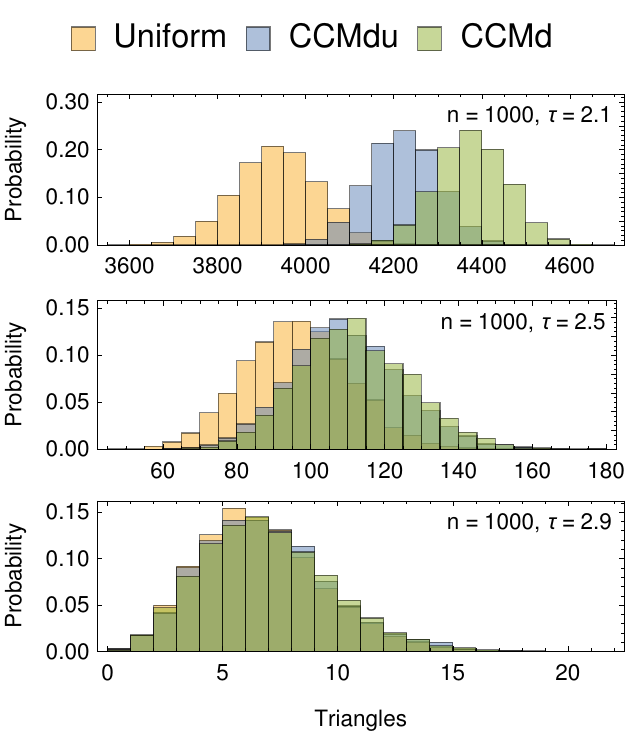}
    \caption{\label{fig:ccminitialtris} Distribution of the number of triangles in the CCMd/CCMdu construction as well as in the uniform distribution. For every $n,\tau$, the canonical degree sequence was used and 5000 samples are shown from each distribution.}
\end{figure}
Figure \ref{fig:ccminitialtris} shows the distribution of the number of triangles in the graphs generated using the constrained configuration model, compared to the uniform distribution. The initial number of triangles in both CCMdu and CCMd is near the uniform average, though slightly higher, whereas the Havel-Hakimi construction generates graphs where the number of triangles is usually many times higher than average. Starting the switch chain process using the CCM construction \emph{may} therefore give a starting point closer to equilibrium in number of triangles.

\subsection{Mixing time}
\begin{figure}
    \includegraphics[width=0.45\textwidth]{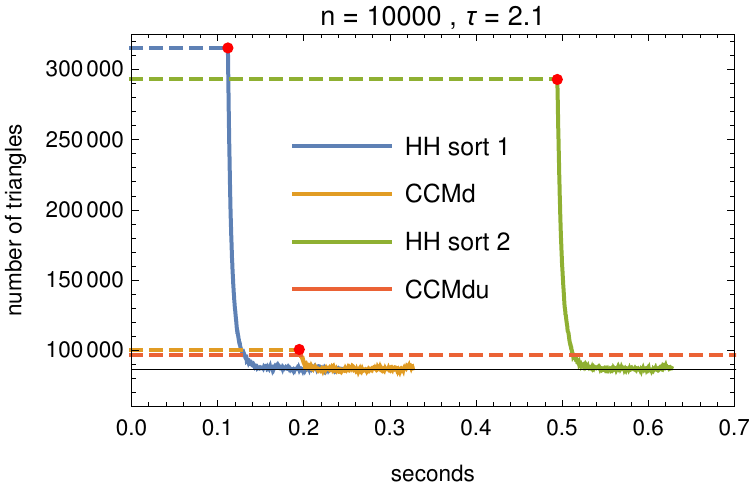}
    \caption{\label{fig:ccm_timeevol} Time evolution of the number of triangles for different initial graphs. The computational time is measured in seconds and the same canonical degree sequence was used for all runs. HH stands for the Havel-Hakimi algorithm and is done twice: using insertion sort (sort 1) and using the C++ standard library sorting algorithm (sort 2). See also section C of the supplementary material. The CCMdu construction took more than 2 full seconds so falls outside the plotrange.}
\end{figure}

Figure \ref{fig:ccm_timeevol} shows the time evolution of the number of triangles for both the Havel-Hakimi starting point and the CCMd and CCMdu starting points, where time is measured in seconds instead of switch chain steps. The plot shows that the construction of the initial graphs takes up a significant portion of the time compared to the switches needed to reach equilibrium. We see that the number of triangles in the Havel-Hakimi starting point lies much further away from the average than that in the CCMd and CCMdu starting points, but still the Havel-Hakimi construction including the mixing time is faster than the CCMd construction, provided a proper sorting algorithm is used. Note that the sorting algorithm influences the number of triangles in the starting graph. The reasons for this are explained in section C of the supplementary material.
The CCM construction (both the CCMd and CCMdu variant) is more computationally intensive and more complicated than the Havel-Hakimi algorithm because it has to keep track of which vertices can be paired to and select a random vertex weighted by the number of remaining half-edges it has. Depending on the implementation (and things like the pseudorandom number generator that is used) it can be faster to start with the simple Havel-Hakimi algorithm and do some extra switches which are much simpler. Note that it might be possible that a faster implementation of the CCMd algorithm beats the Havel-Hakimi algorithm.

\section{Conclusion}
We propose triangle counts as a measure to quantify how close a distribution of simple graphs is to the stationary distribution. Triangles form the simplest non-trivial subgraphs, and contain a large amount of information about the structure of the graph. 
Figure~\ref{fig:timeevol} clearly shows that for scale-free networks with degree exponent $\tau\in(2,3)$, the number of triangles is fluctuating wildly for different degree sequences with the same value of $\tau$.  When two uniform random graphs with the same degree sequences are created, the number of triangles in these two graphs will be close.

The method of choice to simulate a uniform graph with prescribed degrees is the switch chain. When the switch chain is set up properly, its stationary distribution is uniform.
We investigate the role of the starting point of the switch chain. As can perhaps be expected, Havel Hakimi starts from a triangle count that is much higher than for uniform random graphs (even quite close to its maximal value). Instead, we investigated the constrained configuration model CCMd as a starting point. CCMd mimics the configuration model, while ensuring simplicity of the graph, whereas the configuration model could potentially fail to produce a simple graph. Our simulations show that CCMd almost always succeeds, but is computationally heavier and therefore it is faster to use the Havel-Hakimi algorithm and do extra switches. The related CCMdu construction, in which the vertices are ordered by their remaining degrees, is computationally even heavier, while remarkably also having a lower success probability.

Our simulations clearly show that the triangle count for $\tau\in(2,3)$ in the uniform random graph substantially deviates from the often used erased configuration model, where self-loops and multiple edges in the configuration model are removed. We conjecture that the number of triangles in the ECM is higher than the number of triangles in uniform random graphs when the graph size is large enough. Thus, care is needed in the analysis of uniform graphs in the omnipresent scale-free regime when $\tau\in(2,3)$. 

In the mathematical literature, the switch chain is also used as a key methodology to rigorously prove properties of uniform graphs. This technique is limited to cases where at most one triangle is created or destroyed per switch chain move. Our simulations clearly show that this often fails, particularly for $\tau$ small, thus implying that such proofs are doomed to fail.

\acknowledgements
The work of RvdH and TB is supported by the NWO Gravitation Networks grant 024.002.003.
The work of RvdH and CS was supported by NWO TOP grant 613.001.451.
The work of RvdH is further supported by the NWO VICI grant 639.033.806. 

\bibliographystyle{eplbib}
\bibliography{triangle_switch}


\section{Supplementary material}
In this supplementary material we provide more details on the numerical experiments covered in the main manuscript.

\section{A. Aperiodic Markov Chain}
The switch chain on degree sequence $(d_i)_{i\in[n]}$ is aperiodic as long as a simple graph on this degree sequence exists with a path of length 3. When we would choose the outer two edges of such a path for a switch, then this would result in a double edge. Then the Markov Chain would remain in the same state, so that there is a self-loop in the Markov Chain in this state. Therefore, the Markov Chain is aperiodic. Note that when $\Prob{D\geq 2}>0$, this condition is satisfied w.h.p. for $n$ large enough. For more information on the conditions for uniform samples from the switch chain, see \cite{carstens2017}.

\section{B. Sampling from power laws}
We want to sample degree sequences from a distribution with probability mass function 
\begin{equation} \label{eq:degreedistr1}
    \Prob{X=x}=Cx^{-\tau},
\end{equation}
where $C$ is a normalising constant and $x$ is an integer with $x\geq 1$. The constant $C$ is given by $C=1/\zeta(\tau)$ where $\zeta$ is the Riemann zeta function $\zeta(s)=\sum_{n=1}^{\infty} n^{-s}$. To generate samples from such a distribution one would need to compute the constant $C$ and then numerically invert the cumulative distribution function for each sample. This is computationally intensive and hence it is common to make approximations. There are several ways one can do this, the easiest of which is by first sampling from a continuous power-law distribution and then rounding to the nearest integer \cite{clauset2007}, which we did for the simulations presented in the manuscript. The continuous distribution has density proportional to $x^{-\tau}$ for $x\geq 1$ and satisfies $\Prob{X_\mathrm{cont}\leq x} = 1-x^{-(\tau-1)}$. The canonical degree sequence is given by
\begin{equation*}
    d^{(n,\tau)}_i = \left[ \left(\frac{i}{n}\right)^{-\frac{1}{\tau-1}} \right]
\end{equation*}
where $i$ runs from $1$ to $n$ and the rounding $\left[ \cdot \right]$ is to the \emph{nearest} integer. The largest degree in this canonical degree sequence is therefore equal to $\left[n^{\frac{1}{\tau-1}}\right]$.

\section{C. Triangles in Havel-Hakimi construction}
In the random graphs that we consider, the majority of the triangles is due to the high-degree vertices pairing up with each other. By construction, the Havel-Hakimi (Erd\H{o}s-Gallai) algorithm pairs up the high-degree vertices \emph{only} with the other high-degree vertices, and therefore the number of triangles in the resulting graph is larger than the average. In fact, the data suggested that this construction might yield the \emph{maximum possible} number of triangles of all graphs with a given degree sequence. However, Figure 8 in the main manuscript shows an example where two runs of the Havel-Hakimi algorithm (using different sorting algorithms) give a graph with a different number of triangles.
To investigate this, we iterated over all possible graphs with $n$ vertices ($2^{\binom{n}{2}}$ graphs) for $n\in\{5,6,7,8\}$ to compute the maximum number of triangles for every valid degree sequence of size $n$. Then for every valid degree sequence we did the Havel-Hakimi construction to compare the number of triangles. We found that for most degree sequences (1022 out of 1213 valid degree sequences for $n=8$) the Havel-Hakimi construction indeed gave the highest possible number of triangles. However there were some degree sequences for which the construction yielded a graph with fewer triangles than the maximum. This happened for example for the degree sequence $\{4,4,3,3,3,2,1\}$ for which the maximum number of triangles possible is 5. However, in the Havel-Hakimi construction, one first pairs up a degree-4 vertex after which the remaining degrees are $\{0,3,2,2,2,2,1\}$. Then the degree-3 vertex (which was first the degree-4 vertex) is paired up, to three of the four degree-2 vertices. The Havel-Hakimi construction does not specify which three vertices to pick in this case (any ordering will work). Depending on how the vertices are now sorted relative to each other, the construction can result in a graph with 3 triangles instead of 5. This shows that the Havel-Hakimi construction does not always result in the maximum number of triangles if we fix a certain sorting method. However, on the other hand it is still possible that by choosing a specific ordering of vertices, the Havel-Hakimi construction \emph{does} yield the maximum number of triangles possible. For the 191 out of 1213 valid degree sequences where the maximum was not obtained, on average the Havel-Hakimi construction produced $1.57$ fewer triangles than the maximum possible (average only over those 191). Overall we can conclude that the Havel-Hakimi algorithm produces graphs with close-to-maximum number of triangles.

\section{D. Success rate of switch chain moves}
\begin{figure}
    \centering
    \includegraphics[width=0.45\textwidth]{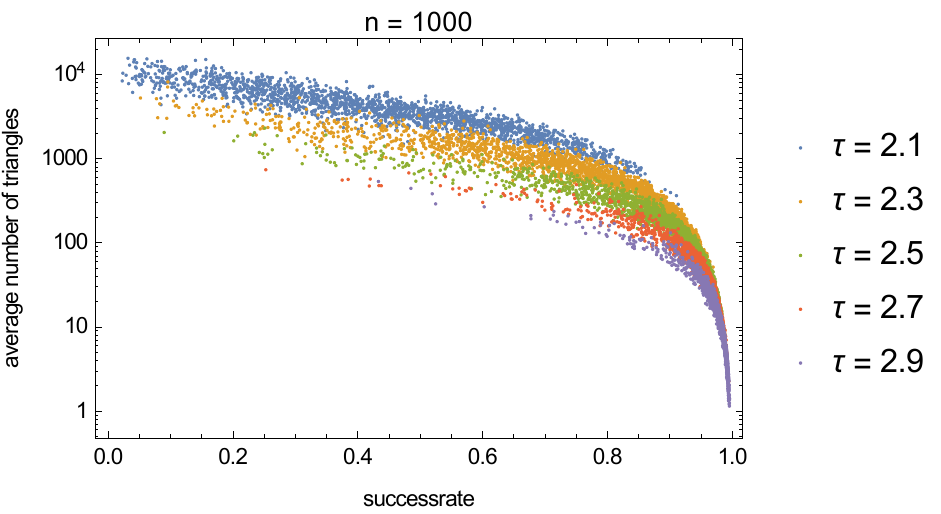}
    \caption{\label{fig:successcorrelation} The plot shows the correlation of the success rate versus the number of triangles in the graph on a logarithmic scale. For every value of $\tau$ shown, $2000$ points are plotted, each corresponding to a degree sequence sampled according to \eqref{eq:degreedistr1} described in section B.}
\end{figure}The switch chain samples two non-touching edges uniformly at random and then tries to `switch' them. This switch step can fail if there were already edges present between the vertices involved in which case the graph is left unchanged (this does count as a step of the Markov Chain, it is a self-loop). We investigated how often this happens, i.e. the success rate of the switch chain. It turns out that the success rate is correlated with the number of triangles present in the graph when performing a switch (not only the number of triangles in equilibrium): the more triangles there are, the lower the success rate. Figure \ref{fig:successcorrelation} shows this correlation.


\section{E. Construction success rates of CCM}

In the main manuscript, we looked at the construction success rates for the CCMd and CCMdu constructions. The conclusion was that the CCMd construction is almost always successful whereas the CCMdu construction is not. Here we show an additional plot to show this.
\begin{figure}
    \centering
    \includegraphics[width=0.45\textwidth]{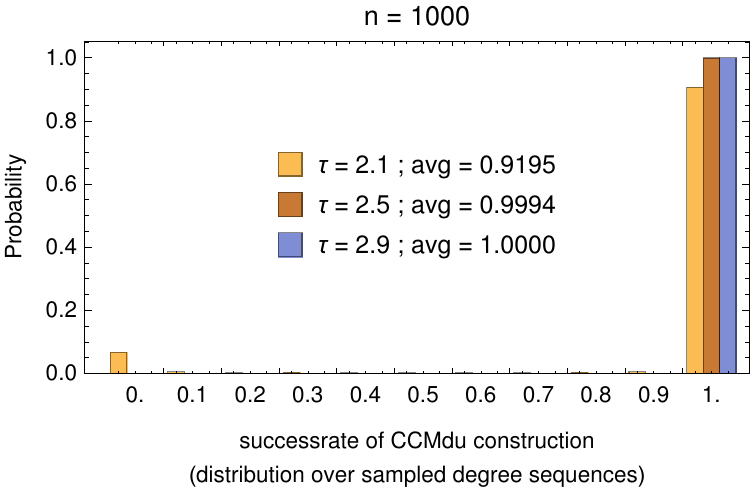}
    \caption{\label{fig:ccmconstructionrate} Construction success rate for CCMdu for $n=1000$ and $\tau\in\{2.1,2.5,2.9\}$. The distribution is over 5000 graphical sampled degree sequences (so not the canonical degree sequence) for each $\tau$, with 200 construction attempts per sequence to determine the successrate. The rightmost column (with label 1.) corresponds to successrates of at least $0.95$.}
\end{figure}
Figure \ref{fig:ccmconstructionrate} shows the construction success rates for the CCMdu construction. Interestingly, for $\tau$ close to 2, a degree sequence is either very `good' or very `bad'. The success rate is is less than $0.05$ for some degree sequences (meaning less than 10 out of 200 attempts succeeded) and higher than $0.95$ for most. The figure also shows that the success rate increases with $\tau$.

\end{document}